\documentclass[reqno]{elsart}
\usepackage{graphicx,amssymb,amsmath,subfigure}

\newenvironment{defi}{\noindent {\bf Definition}}{}
\newenvironment{proof}{\emph{Proof:\quad}}{\qed\bigskip}

\newcommand{\T}{\mathcal{T}}
\newcommand{\PT}{\mathcal{PT}}
\newcommand{\PPT}{\mathcal{PPT}}

\graphicspath{{figures/}}

\begin{document}

\begin{frontmatter}

\title
{On the Number of Pseudo-Triangulations \\ of Certain Point Sets
\thanksref{title}
}

\date{14 June 2007}

\thanks[title]{Parts of this work were done while the authors visited the
Departament de Matem\`{a}tica Aplicada~II, Universitat
Polit\`{e}cnica de Catalunya, the Institut f\"ur
Softwaretechnologie, Technische Universit\"at Graz and the
Departamento de Matem\'a\-ticas, Universidad de Alcal\'a. \\
Research
of Oswin Aichholzer is partially supported by Acciones Integradas
2003-2004, Proj.Nr.1/2003 and by the FWF [Austrian Fonds zur
F\"orderung der Wissenschaftlichen Forschung] under grant
S09205-N12, FSP Industrial Geometry. Research of David Orden is
partially supported by grants HU2002-0010, MTM2005-08618-C02-02 and
S-0505/DPI/0235-02. Research of Francisco Santos is partially
supported by Acci\'on Integrada Espa\~na-Austria HU2002-0010 and
grant MTM2005-08618-C02-02 of Spanish Direcci\'on General de
Investigaci\'on Cient\'{\i}fica.
}

\author{Oswin Aichholzer}
\address{Institute for Software Technology, Graz University of Technology, Austria.} 
\ead{oaich@ist.tugraz.at}

\author{David Orden}
\address{Departamento de  Matem\'aticas,  Universidad  de  Alcal\'a, Spain.}
 \ead { david.orden@uah.es}

  \author{Francisco Santos}
\address{Depto. de Matem\'aticas, Estad\'{\i}stica y Computaci\'on, Univ. de Cantabria, Spain.}
 \ead{ santosf@unican.es}

 \author{Bettina Speckmann}
\address{Dept. of Mathematics and Computer Science, TU Eindhoven, The Netherlands.}
\ead{\tt speckman@win.tue.nl}



\begin{abstract}
We pose a monotonicity conjecture on the number of pseudo-triangulations of any planar
point set, and check it on two prominent families of point sets, namely the so-called
double circle and double chain. The latter has asymptotically $12^n n^{\Theta(1)}$
pointed pseudo-triangulations, which lies significantly above the maximum number of
triangulations in a planar point set known so far.
\end{abstract}


\end{frontmatter}
\section{Introduction}

Pseudo-triangulations, also called geodesic triangulations, are a
generalization of triangulations which has found multiple applications in Computational
Geometry in the last ten years. They were originally introduced in the context of
visibility~\cite{pv-mtvg-96,pv-tsvcpt-96} and ray
shooting~\cite{cegghss-rspug-94,gt-drssp-97}, but recently have also been applied in
kinetic collision detection~\cite{abghz-dfstkcd-2000,kss_kcdsp_01}, and guarding
problems~\cite{st03:avsppt}, among others. They also have surprising relations to
rigidity~\cite{horsssssw,osss-cpt-05,streinu-00,streinu-05} and locally convex
functions~\cite{aabh03}. See the recent survey~\cite{RSS07} for more information.

\begin{figure}[htb]
  \centering
  \includegraphics{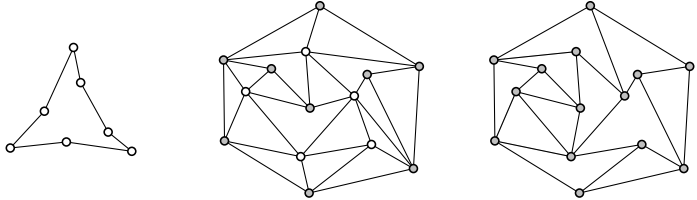}
  \caption{A pseudo-triangle (left), a pseudo-triangulation (middle), a pointed pseudo-triangulation (right).}
  \label{fig:pseudo}
\end{figure}

\begin{figure}
 \centering
 \includegraphics{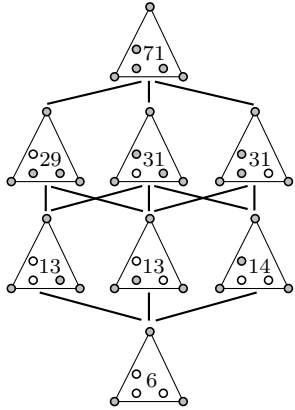}
 \caption{Hasse diagram of the subsets $W$ of $A_{I}$ for a set of 6 points.}
 \label{fig:hasse}
\end{figure}

A {\em pseudo-triangle} is a planar polygon that has exactly three
convex vertices with internal angles less than $\pi$. These vertices
are called \emph{corners} and the three inward convex polygonal
chains joining them are called \emph{pseudo-edges} of the
pseudo-triangle. A {\em pseudo-triangulation} for a set $A$ of $n$
points in the plane is a partition of $\operatorname{conv}(A)$ into
pseudo-triangles whose vertex set is exactly~$A$. Although
pseudo-triangulations can be studied for general point
sets~\cite{OS}, in this paper we will consider only point sets in
general position. A vertex is {\em pointed} if it has an incident
angle greater than~$\pi$. A {\em pointed pseudo-triangulation} is a
pseudo-triangulation where every vertex is pointed. See, for
example, Figure~\ref{fig:pseudo}---here, and in
Figure~\ref{fig:hasse}, pointed vertices are dark, non-pointed
vertices are light. Note that the vertices of
$\operatorname{conv}(A)$ are always pointed.

The set of all pseudo-triangulations of a point set has somewhat nicer properties than
that of all triangulations. For example, pseudo-triangulations of a point set with $n$
elements form the vertex set of a certain polytope of dimension $3n-3$ whose edges
correspond to flips \cite{OS}. The diameter of the graph of pseudo-triangulations is
$O(n\log n)$ \cite{aak03} versus the $\Theta(n^2)$ diameter of the graph of
triangulations of certain point sets. Also, for standard triangulations, it is not known
which sets with a given number of points have the fewest or the most triangulations, but
it was shown in~\cite{aaks-cmpt-02} that sets in convex position minimize the number of
pointed pseudo-triangulations among all point sets with a given number of vertices (hence the number of all pseudo-triangulations, since in
convex position all pseudo-triangulations are pointed).

\medskip
 Let $A$ be a point set and let $A_I$ be the subset of its interior points. Let
$\PT(A)$ be the set of pseudo-triangulations of $A$. This set can naturally be stratified
into $2^{A_I}$ sets, one for each possible subset of $A_I$. More precisely, for each
subset $W\subseteq A_I$ we denote by $\PT_{W}(A)$ the set of pseudo-triangulations of $A$
in which the points of $W$ are pointed and those of $A_I\setminus W$ are non-pointed. For
example, $\PT_{\emptyset}(A)$ is the set of triangulations of $A$, which we abbreviate
as~$\T(A)$. Similarly, $\PT_{A_I}(A)$ is the set of pointed pseudo-triangulations of $A$,
that we abbreviate as $\PPT(A)$.  In \cite{RRSS01}, the following inequality is proved:
for every $W\subseteq A_I$ and every $p\in W$,
\begin{equation}
3 \, |\PT_{W\backslash \{p\}}(A)| \ge |\PT_{W}(A)|. \label{eqn:RRSS01}
\end{equation}
The main goal of this paper is to explore the relation between the numbers of
triangulations, pointed pseudo-triangulations, and everything in between, for several
specific point sets. In particular, we test the following conjecture, which is implicit
in previous work, but stated here explicitly for the first time:
\begin{conj}
\label{conj:monotone} For every point set $A$ in general position in the plane, the
cardinalities of $ \PT_{W}(A)$ are monotone with respect to $W$. That is to say, for any
subset $W$ of $A$'s interior points and for every $p\in W$, one has
\[
| \PT_{W}(A)| \ge | \PT_{W\backslash\{p\}}(A)|.
\]
\end{conj}

The conjecture is consistent with the following result, also from \cite{RRSS01}: If $A$
has a single interior point and $n-1$ boundary points, then $|\PPT(A)|$ is  greater than
$|\T(A)|$. Actually, the difference is always equal to the Catalan number $C_{n-2}$, no
matter where the interior point is, while $|\T(A)|$ ranges from $C_{n-2}-C_{n-3}\simeq
\frac{3}{4}C_{n-2}$ when the interior point is near the boundary to essentially $C_{n-2}$
when it is near the center.

Conjecture~\ref{conj:monotone} does not imply that the number of
pseudo-triangu\-la\-tions of $A$ with, say, $k$ pointed vertices is greater than the
number of them with $k-1$ pointed vertices. For example, Figure~\ref{fig:hasse} shows
the eight possibilities of $| \PT_{W}(A)|$ for a set of six points, three of them
interior, displayed in the Hasse diagram of subsets of $A_I$. The numbers satisfy
Conjecture~\ref{conj:monotone}, but there are less pointed pseudo-triangulations (71)
than pseudo-triangulations with one non-pointed and two pointed vertices ($29+31+31=91$).

Actually, formula~(\ref{eqn:RRSS01}) says that the same will happen
for any point set with at least four interior points. Applied with
$W=A_I$ and taking the different possibilities for $v\in A_I$, the
formula gives
\[
\sum_{v\in A_I} 3 \, |\PT_{A_I\backslash \{v\}}(A)| \ge  |A_I|\cdot |\PT_{A_I}(A)|.
\]
In other words, the ratio of pseudo-triangulations with exactly one non-pointed vertex to
pointed pseudo-triangulations is at least $|A_I|/3$.

Similarly, the monotonicity is conjectured  with respect to the sets $W$ and not only
their cardinalities: There exists a set $A$ of 10 points, 7 of them interior, and two
subsets $W$ and $W'$ of four and three interior points, respectively, with $| \PT_{W}(A)|
<| \PT_{W'}(A)|$. There are no examples like this with less than 10 points.

As initial evidence for Conjecture~\ref{conj:monotone} we have
computed the numbers of triangulations and of pointed
pseudo-triangulations for all order types of planar point sets in
general position with 10 points or less. This has been done using
the order type database in~\cite{AK01}. One implication of
Conjecture~\ref{conj:monotone} is that every point set should have
at least as many pointed pseudo-triangulations as triangulations,
and this is actually the case up to 10 points. Even more, the
natural expectation is that the ratio between those two numbers
grows exponentially with the number of interior points~$i$, the base
of the exponent being between 1 (by Conjecture~\ref{conj:monotone})
and 3 (by inequality~(\ref{eqn:RRSS01})). In Tables \ref{table:max}
and \ref{table:min} we show the maximum and the minimum values of
the ratio $(|\PPT(A)|/|\T(A)|)^{1/i}$ obtained for each value of the
total number of points $n$ and of interior points~$i$.
\begin{table}[htb]
\begin{center}
\begin{tabular}{|c|ccccccc|}
\hline
   & $i=1$   &       $i=2$    &     $i= 3$     &    $i= 4$    &      $i=5$     &    $i= 6$   &   $ i= 7$ \cr
\hline $n-i=3$   &  3.00000 &  2.54951 &  2.50665 &  2.38010 &  2.31659 &  2.26583  &
2.23025 \cr $n-i=4$   &  2.66667 &  2.51661 &  2.47042 &  2.44151 &  2.35995 &  2.30562 &
\cr $n-i=5$   &  2.55556 &  2.48151 &  2.44824 &  2.42734 &  2.41308 &   & \cr $n-i=6$
&  2.50000 &  2.45607 &  2.43210 &  2.41625 &   & & \cr $n-i=7$   &  2.46667 &  2.43763 &
2.41980 &   & & & \cr $n-i=8$   &  2.44444 &  2.42384 &   & & & & \cr $n-i=9$   &
2.42857 &   & & & & & \cr \hline
\end{tabular}
\caption {Maximum values of $\left(\frac{|\PPT(A)|}{|\T(A)|}\right)^{1/i}$ for order
types with at most 10 points.} \label{table:max}
\end{center}
\end{table}
\begin{table}[htb]
\begin{center}

\smallskip
\begin{tabular}{|c|ccccccc|}
\hline
          &  $i=1$   &      $i= 2$    &    $ i= 3$     &   $ i= 4$    &      $i=5$     &    $i= 6$  &    $i= 7$  \cr
\hline $n-i=3$  &  3.00000 &  2.54951 &  2.11791  & 2.00415 &  1.88343  & 1.80952  &
1.75590\cr $n-i=4$  & 2.66667  &  2.29129 &  2.01550  & 1.91670 &  1.82364  & 1.76240
&\cr $n-i=5$  &  2.27273 &  2.08637 &  1.91798  & 1.84371 &  1.77002 & &\cr $n-i=6$  &
2.16667 &  1.99211 &  1.83797 &  1.78048 & & &\cr $n-i=7$  &  2.03937 &  1.91361 &
1.80419 & & & &\cr $n-i=8$  &  1.98621 &  1.86445 & & & & &\cr $n-i=9$  &  1.92318 &  & &
& & & \cr \hline
\end{tabular}
\caption {Minimum values of $\left(\frac{|\PPT(A)|}{|\T(A)|}\right)^{1/i}$ for order
types with at most 10 points.} \label{table:min}
\end{center}
\end{table}
It is interesting to observe that rows (fixed number of boundary points) and columns
(fixed number of interior points) are monotone in both tables, while diagonals (fixed
total number of points) are not always monotone in Table~\ref{table:min} (see diagonal
$n=6$ and $n=9$).

From Table \ref{table:asymptotic} below we can derive the (asymptotic) value of the
same parameter $\left({|\PPT(A)|}/{|\T(A)|}\right)^{1/i}$ for certain families of planar
point sets which are the main object in this paper: double circle, single chain, and
double chain. The results are $7/3 \approx 2.333$, $2$ and $1.5$, respectively. Also, the results
in~\cite{RRSS01} say that if $A$ has a single interior point then
\[
\frac{7}{3}\simeq 1+\frac{C_{n-2}}{C_{n-2}-C_{n-3}} \ge {|\PPT(A)|}/{|\T(A)|} \ge 1
+\frac{ 2C_{n-2}}{C_{n-1}-(n-1)C_{(n-2)/2}^{{}^2}}\simeq 2,
\]
with equality on the left  when the interior point is close to the boundary, and on the
right when the interior point is at the center of a regular $(n-1)$-gon, with $n$ even.

\bigskip 

\begin{figure}[htb]
  \centering
  \includegraphics{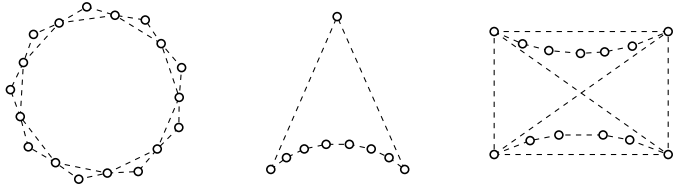}
  \caption{A double circle (left), a single chain (middle), a double chain (right).}
  \label{fig:pointsets}
\end{figure}

 In the rest of the paper we  consider three
families of point sets in the plane: ``double circles'', ``double
chains'', and what we call ``single chains''. See
Figure~\ref{fig:pointsets} for examples, the exact definitions are
given in the respective sections. The double circle is conjectured
to be the point set with asymptotically the smallest number of
triangulations, for a fixed number of points~\cite{AK01}. The double
chain has been the example with (asymptotically)  the biggest number
of triangulations known~($\Theta^*(8^n)$, see~\cite{SS03}), until a
new structure was found recently, the so-called double zig-zag
chain, with $\Theta^*(\sqrt{72}^n)$ triangulations~\cite{ahhhkv06}.
(Here and in the rest of the paper, the notation $\Theta^*$ means
that a polynomial factor is neglected.) We studied single chains
originally as a step to analyze double chains, but it turns out,
that the number of pseudo-triangulations of single chains also has
very interesting combinatorial properties (see
Theorem~\ref{thm:singlechain}).

Our interest in the number of pseudo-triangulations of these point
sets is two-fold. On the one hand, we prove that
Conjecture~\ref{conj:monotone} holds in these three cases. On the
other hand, we are interested in how many pseudo-triangulations a
general point set in the plane can have. The study of these point
sets, which have very many or very few triangulations, should give
an indication of it. Even if the minimum number of
pseudo-triangulations is achieved by the convex $n$-gon (as
mentioned above), it may well be that the double circle (or, more
generally, the point sets in ``almost convex position'' studied in
Section~\ref{sec:doublecircle}) minimize the numbers of
pseudo-triangulations for fixed numbers of boundary and interior
points.

Our main results are summarized in Table~\ref{table:asymptotic}, where we show only the
global number of all pseudo-triangulations and the extremal cases of triangulations and
pointed pseudo-triangulations. In all cases $n$ is assumed to be the total number of
points, and a factor polynomial in $n$ has been neglected.  The double circle has $n/2$
interior points and the single and double chains have $n-3$ and $n-4$ interior points,
respectively.
\begin{table}[htb]
\begin{center}
\begin{tabular}{|c||c|c|c|}
\hline
   & double  & single  & double  \\
   & circle    & chain   & chain   \\
\hline $|\T(A)|$
    &$\sqrt{12}^n$ & $4^n$    &   $8^n$\\
\hline $|\PPT(A)|$
    &$\sqrt{28}^n$   &   $8^n$ &    $12^n$\\
\hline $|\PT(A)|$
    &$\sqrt{40}^n$   &   $12^n$ &   $20^n$\\
\hline Conjecture \ref{conj:monotone}
      & Holds       &        Holds      &     Holds  \\
\hline
\end{tabular}
\medskip
\end{center}
\caption{ Asymptotic number of triangulations, pointed pseudo-triangulations and
pseudo-triangulations, for special point sets.} \label{table:asymptotic}
\end{table}

The paper is organized as follows: Section~\ref{sec:doublecircle}
studies the number of pseudo-triangulations of so-called point sets
in almost convex position, among which the double circle is the
extremal case. The next three sections are devoted to the single
chain. Section~\ref{sec:singlechain} gives approximations, within a
factor of four, for the numbers of pseudo-triangulations $|\PT(A)|$
and pointed pseudo-triangulations $|\PPT(A)|$ of the single chain.
The proof of the crucial result that gives the asymptotics,
Theorem~\ref{thm:(l+3)-gon}, is given separately in
Section~\ref{sec:conjecture3}. Section~\ref{section:extras} uses a
different approach to provide a much better approximation
of~$|\PT(A)|$ and~$|\PPT(A)|$ for the single chain. Finally,
Section~\ref{sec:doublechain} is devoted to the double chain, whose
study is based on that of the single chain and, in particular, on
the aforementioned Theorem~\ref{thm:(l+3)-gon}.

\section{The double circle and its relatives}
\label{sec:doublecircle}

 For any given pair of positive integers $v\ge 3$ and $i\le
v$, we say that a point set $A$ is {\em in  almost convex position with
parameters $(v, i)$} if it consists of a set of $v$ points forming
the vertex set of a convex $v$-gon and a set of $i$ interior points,
placed ``sufficiently close'' to $i$ different edges of the $v$-gon.
Here, we say that an interior point $p$ is placed \emph{sufficiently close to the edge
$(r, q)$ of the $v$-gon} if no segment connecting two points of $A$
 can separate $p$ from $(r, q)$. The {\em double circle} is
the extremal case with $v=i=n/2$, where there is one interior point
close to every boundary edge. It has asymptotically
$\Theta(\sqrt{12}^n n^{-3/2})$ triangulations \cite{SS03} and it is
conjectured in \cite{AK01} that this is the smallest number of
triangulations that $n$ points in general position in the plane can
have. This conjecture is known to be true for $n \leq
11$~\cite{ahn04}.

Point sets in almost convex position are a special case of what is
called ``almost-convex polygons'' in~\cite{HN97}. There it is shown
that the number of triangulations of such a point set does not
depend on the choice of the $i$ edges of the $v$-gon.  Indeed, if we
call this number $t(v,i)$, the case $W=\emptyset$ of  Lemma~\ref{lemma:doublecircle}
below provides the recursive formula
\begin{equation}
  t(v,i) = t(v+1,i-1)-t(v,i-1)
\label{eqn:doublecircletriangs}
\end{equation}
which allows to compute $t(v,i)$ starting with $t(v,0)=C_{v-2}$
(Catalan numbers).
It is interesting that  formula~(\ref{eqn:doublecircletriangs}) can be applied to
generate $t(v,i)$ even for $i>v$. The
array obtained by this recursion (difference array of Catalan
numbers) appears in Sloane's Online Encyclopedia of Integer
Sequences \cite{integersequences} with ID number A059346.
The numbers obtained for $i>v$ do not have a
meaning as triangulations of point sets, but (for small values of
$v$) they have other combinatorial interpretations. For example, the
sequence $M_n:=t(3,n)$ forms the Motzkin numbers (number of lattice
paths from $(0,0)$ to $(n,0)$ with steps $(1,0)$, $(1,1)$ or
$(1,-1)$ and lying above the horizontal axis,
see~\cite{ds-77,stanley-ec2}). In the proof of Corollary~\ref{coro:doublecircle} below we use
their asymptotic expression, which appears
for example in~\cite[Section VI.4]{fs-07}:
\begin{equation}
\label{eqn:motzkin}
M_n =  \sqrt{\frac{3}{4\pi}} 3^n \left(n^{-3/2} +O(n^{-5/2})\right) \in \Theta(3^n n^{-3/2}).
\end{equation}

\begin{figure}[htb]
 \centering
 \includegraphics{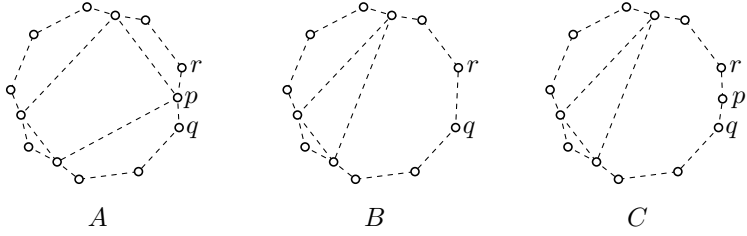}
 \caption{Almost convex point sets: set $A$ with $v = 9$ and $i = 4$, set $B$ with $v = 9$ and $i = 3$, and set $C$ with $v = 10$ and $i = 3$.}
 \label{fig:almconvABC}
\end{figure}

We now generalize the recursive formula (\ref{eqn:doublecircletriangs}) to deal also with
pseudo-triangu\-la\-tions. Let $p$ be a specific interior point of a point set $A$ in
almost convex position and let $(q, r)$ be the convex hull edge which has $p$ next to it. Let $B$ and $C$ be the point sets obtained respectively by deleting $p$ from $A$ and by
moving $p$ to convex position across the edge $(q, r)$ (see Figure~\ref{fig:almconvABC}).

\begin{lem}
\label{lemma:doublecircle} For every $W\subseteq A_I$ not
containing $p$ (so that $W$ is also a set of interior points of $B$ and $C$) one has:
\begin{enumerate}
\item $|\PT_W(A)| = |\PT_W(C)| - |\PT_W(B)|$.
\item $|\PT_{W\cup\{p\}}(A)| = 2 \,|\PT_W(C)| - |\PT_W(B)|$.
\item $|\PT_{W\cup\{p\}}(A)| =  2 |\PT_{W}(A)| + |\PT_{W}(B)|$. In particular, $A$ satisfies Conjecture
\ref{conj:monotone}.
\end{enumerate}
\end{lem}

\begin{proof}
  It is clear that there are bijections between: $(i)$
  pseudo-triangulations of $C$ pointed at $W$ that use the edge $(q, r)$
  and pseudo-triangulations of $B$ pointed at $W$ and $(ii)$
  pseudo-triangulations of $C$ pointed at $W$ that do not use the edge
  $(q, r)$ and pseudo-triangulations of $A$ pointed at $W$.  These
  bijections prove part (1).

  To prove part (2), we partition the pseudo-triangulations of $A$ in which $p$ is pointed into three sets: those using the edges $(p, q)$ and $(p, r)$ (and hence having
  no other edge incident to $p$), those using $(p, q)$ but not $(p, r)$, and
  those using $(p, r)$ but not $(p, q)$. The first set is in bijection with
  the pseudo-triangulations of $B$.  Each of the other two is in
  bijection with pseudo-triangulations of $C$ that do not use the edge
  $(q, r)$, that is pseudo-triangulations of $C$ minus those of $B$.
  Since the bijections preserve pointedness at interior points (other
  than $p$), we get $|\PT_{W\cup\{p\}}(A)| = |\PT_W(B)| + 2 \,(|
  \PT_W(C)| - |\PT_W(B)|)$, as desired.

  Part (3) is obtained eliminating $|\PT_W(C)|$ from parts (1) and (2).
\end{proof}

This lemma shows that $|\PT_{W}(A)|$ only depends on the parameters $(v,i)$ of $A$ and the number $|W|$ of points prescribed to be interior. Indeed,
let us call $s(v,j,k)=|\PT_{W}(A)|$ where $v$ is the number of boundary points and $k=|W|$ and $j=i-k$ are the numbers of interior points prescribed to be pointed and non-pointed, respectively. Then,  parts (1) and  (2) of the lemma translate to
\begin{equation}
s(v,j+1,k) = s(v+1, j,k) - s(v , j,k),
\label{eq:dcircle-1}
\end{equation}
and
\begin{equation}
s(v,j,k+1) = 2 s(v+1, j,k) - s(v , j , k).
\label{eq:dcircle-2}
\end{equation}

From this, all the numbers can be computed recursively: the second formula allows to compute them from the ones with $k=0$ and the first formula allows to compute those from the numbers $s(v,0,0)=C_{v-2}$. Also, from the fact that $C_{v-2}\in \Theta^*(4^v)$ the formulas lead easily to the guess
that $s(v,j,k)\in \Theta^*(4^v 3^j 7^{k})$, which we now prove:

\begin{cor}
\label{coro:doublecircle} The number $s(v,j,k)=|\PT_W(A)|$, where $A$ is a point set in almost convex position with $v$ boundary and $j+k$ interior points, and $W$ is a subset of $k$ of them, satisfies
\[
\Omega( (v+j+k)^{-3/2} )  \le \frac{s(v,j,k)}{4^v 3^j 7^{k}} \le O( j^{-3/2} ).
\]
\end{cor}

\begin{proof}
We start with the following slightly nicer versions of the recursions (\ref{eq:dcircle-1})
and (\ref{eq:dcircle-2}). The first one is a simple rewrite of (\ref{eq:dcircle-1}) and the second is obtained eliminating $s(v , j , k)$:
\begin{eqnarray}
s(v+1,j,k) &=& s(v, j+1,k) + s(v , j,k),
\\
s(v, j,k+1)   &=&  s(v+1,j,k) +  s(v , j+1 , k).
\end{eqnarray}

If we now define $r(v,j,k) =  \frac{s(v,j,k)}{4^v 3^j 7^{k}}$, these two recursions translate to:
\begin{eqnarray}
r(v+1,j,k) &=& \frac{3}{4}r(v, j+1,k) +  \frac{1}{4}r(v , j,k),
\\
r(v, j,k+1)   &=&  \frac{4}{7}r(v+1,j,k) +   \frac{3}{7}r(v , j+1 , k).
\end{eqnarray}
From this eventually we get
\[
\min_{n=j,\dots, v + j + k -3} r(3,n,0) \le r(v , j,k) \le \max_{n=j,\dots, v + j + k -3} r(3,n,0).
\]
That is to say,
\[
\min_{n=j,\dots, v + j + k -3} \frac{M_{n}}{4^3 3^n} \le \frac{s(v , j,k)}{4^v 3^j 7^{k}} \le \max_{n=j,\dots, v + j + k -3} \frac{M_{n}}{4^3 3^n},
\]
where the sequence $M_n=s(3,n,0)$ are the afore mentioned Motzkin numbers.
The fact that $M_n\in \Theta(3^n n^{-3/2})$ finishes the proof.
\end{proof}

We believe that a finer use of the asymptotics of the Motzkin numbers would lead to the slightly stronger statement that $ s(v,j,k)\in \Theta( 4^v 3^j 7^{k}  (v+j+k)^{-3/2} )$.
Anyway, Corollary~\ref{coro:doublecircle} implies that
\[
|\PPT(A)|\in \Theta^*(4^v 7^i), \qquad \mbox{and} \qquad |\PT(A)|\in \Theta^*(4^v 10^i),
\]
where the last formula comes from adding $|\PT_W(A)|=s(v,j,k)$ over
all the $2^i$ values of $W$:
\[
\sum_{k=0}^i {i \choose k} 4^v 3^{i-k} 7^{k}= 4^v 10^{i}.
\]
In conclusion, the double circle ($i=v=n/2$) has about  $\sqrt{28}^n$ pointed
pseudo-triangulations and $\sqrt{40}^n$ pseudo-triangulations in total, modulo a
polynomial factor.

We close this section deriving direct recurrences for the total numbers of pseudo-triangulations and of pointed pseudo-triangulations of point sets in almost convex position.
\begin{cor}
\label{cor:almostconvexposition}
Let $pt(v,i)$ and $ppt(v,i)$ denote  the numbers of pseudo-triangu\-la\-tions and
pointed pseudo-triangula\-tions of a point set in almost convex position with parameters
$(v,i)$, respectively. Then:
\begin{enumerate}
\item $ppt(v,i) = 2\, ppt(v+1,i-1) -ppt(v,i-1)$.
\item $pt(v,i) = 3\, pt(v+1,i-1) -2\,pt(v,i-1)$.
\end{enumerate}
\end{cor}

\begin{proof}
Part (1) is the case $W\cup\{p\}=A_I$ of part (2) of Lemma~\ref{lemma:doublecircle}.
For part (2) we  add parts (1) and (2) of the same lemma over all possible values of $W$.
\end{proof}


\section{The single chain}
\label{sec:singlechain}

\begin{figure}[htb]
  \centering
 \includegraphics{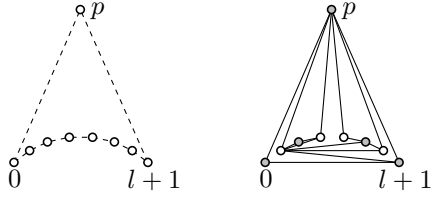}
  \caption{A single chain $A$ with $l=6$ (left), a pointed pseudo-triangulation in~$\PPT_{\{1,3,4,6\}}(A)$ (right).
  }
  \label{fig:PPT_W}
\end{figure}

 Throughout this section let $A$ be the point set with the following $l+3$
points: $l+2$ points labelled $0,1,\dots,l,l+1$ forming a convex $(l+2)$-gon, plus a
vertex $p$ exterior to this polygon and  seeing all edges of it except the edge
$(0,l+1)$. $A$ has three convex hull vertices $\{p,0,l+1\}$ and $l$ interior vertices
$\{1,\dots,l\}$. We call $A$ a {\em single chain} and call $p$ the {\em tip} of $A$, see
Figure~\ref{fig:PPT_W} (left).

Besides classifying pseudo-triangulations of $A$ with respect to their sets of pointed
vertices, here we need to classify the {\em pointed} pseudo-triangulations of the single
chain according to which interior points are joined to the tip. That is, for each subset
$W\subseteq A_I$ we denote by $\PPT_W(A)$ the set of pointed pseudo-triangulations of $A$
in which $p$ is joined to $i$ if and only if $i\in W$, see Figure~\ref{fig:PPT_W} (right).

It is easy to realize that from the numbers $|\PPT_W(A)|$ one can
recover the numbers $|\PT_W(A)|$, which are our main interest in
this section:
\begin{lem}
\label{lemma:ppt} For every $W$:
\[
|\PT_W(A)| =\sum_{W'\subseteq W} |\PPT_{W'}(A)|.
\]
In particular, Conjecture \ref{conj:monotone} holds for the single chain.
\end{lem}
\begin{proof}
  In every pseudo-triangulation of $\PT_W(A)$ the non-pointed vertices $i\in A_I\backslash W$
  have to be joined to the tip~$p$. If we delete those edges $(p,p_i)$ for
  all $i\in A_I\backslash W$, we get an element of a certain
  $\PPT_{W'}(A)$ with $W'\subseteq W$ (here, $W'$ are the vertices
  of the pseudo-triangulation which are joined to $p$ but are pointed).
  This process can clearly be reversed.
\end{proof}

\begin{figure}[htb]
  \centering
 \includegraphics[scale=0.33]{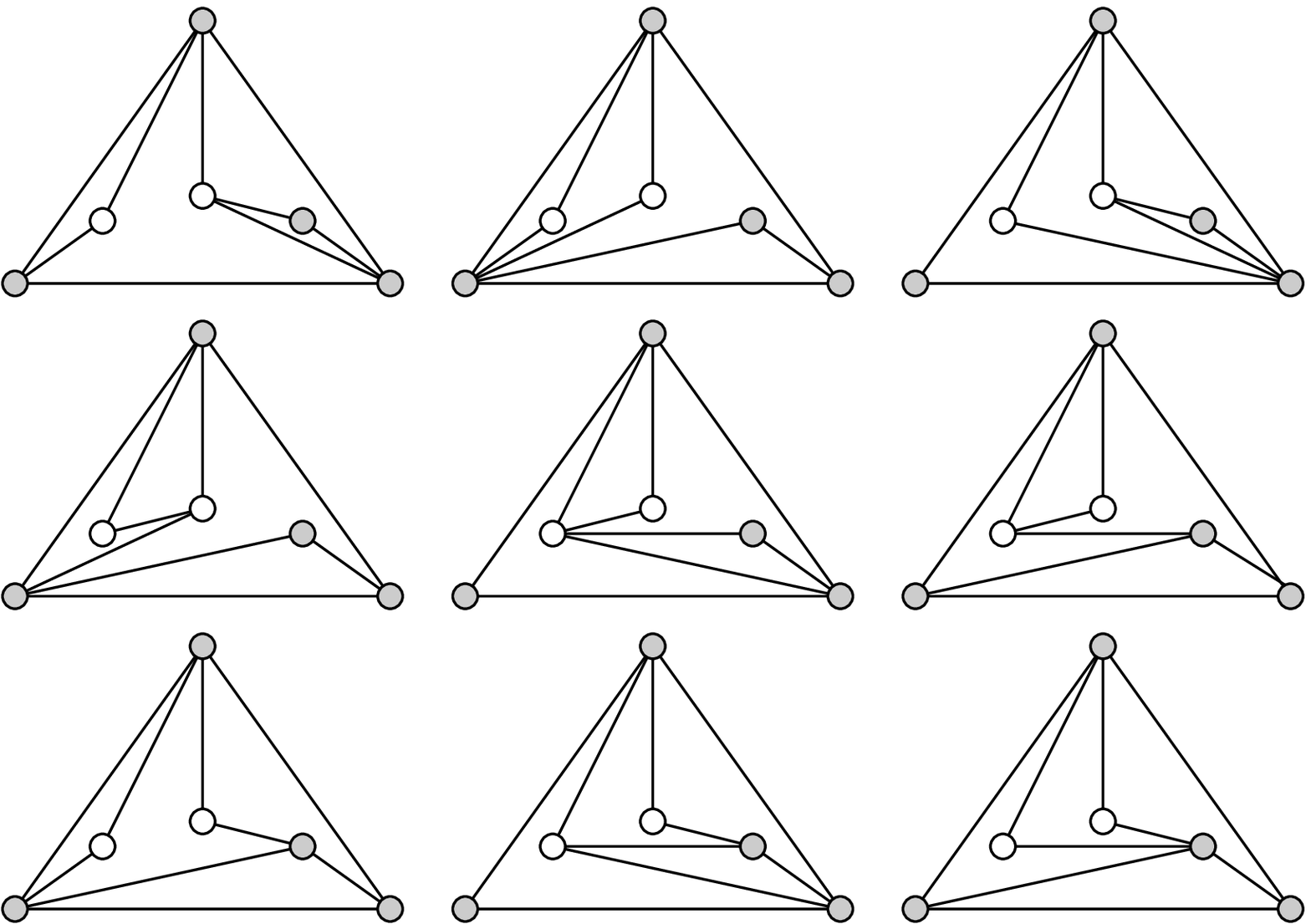}
\hskip 0.5cm
 \includegraphics[scale=0.33]{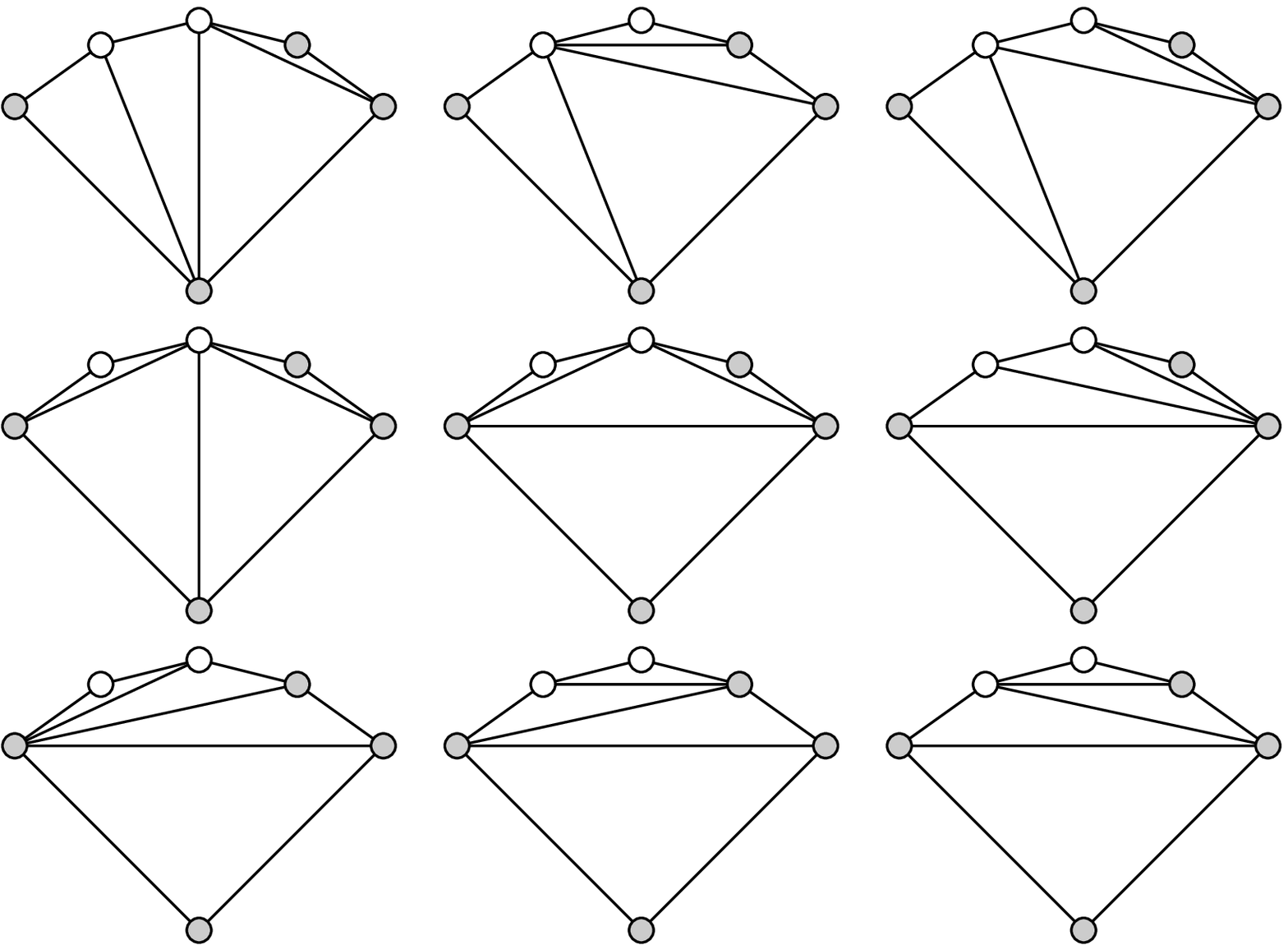}
  \caption{ The nine pointed pseudo-triangulations in $\PPT_W(A)$
  and the nine triangulations in which the interior neighbors of $q$ are contained in $W$.
  In this example $l=3$ and $W=\{1,2\}$.}
  \label{fig:bijection}
\end{figure}

 The following theorem is probably the most surprising
result in this paper. It says that the sets $\PPT_W(A)$ have the
same cardinality as certain subsets of triangulations of a convex
$(l+3)$-gon. Apart of its intrinsic interest, this result
automatically gives the asymptotics of all the numbers $|\PPT_W(A)|$
(Corollary~\ref{coro:singlechain-ppt}), and hence of all the
$|\PT_W(A)|$, too (Corollary~\ref{coro:singlechain-pt}).
\begin{thm}\label{thm:(l+3)-gon} Let $A$ be a single
chain with $l$ interior vertices, let $W\subseteq A_I$ be a subset of them.
 Let $B$ be the convex $(l+3)$-polygon with vertex set
$A\setminus\{p\} \cup \{q\}$, where $q$ is an extra point on the
side opposite to $p$.

The pointed pseudo-triangulations of $A$ in which the interior
neighbors of $p$ are {\bf exactly} the points in $W$ (that is, the
elements of $\PPT_W(A)$) have the same cardinality as the
triangulations of $B$ in which the interior neighbors of $q$ are
{\bf contained} in $W$.
\end{thm}
See Figure~\ref{fig:bijection} for an example.
To maintain the flow of ideas, we postpone the proof of
Theorem~\ref{thm:(l+3)-gon} to Section~\ref{sec:conjecture3}. Let us remark only that our
proof is rather indirect. In particular, it is far from being an explicit bijection
between the two sets involved.

\begin{cor}
\label{coro:singlechain-ppt}
\begin{enumerate}
\item The numbers $|\PPT_{W}(A)|$ are strictly monotone with respect to $W$. That is to
say, for every $W\subset\{1,\dots,l\}$ and every $i\in \{1,\dots,l\} \setminus W$,
\[
| \PPT_{W\cup\{i\}}(A)| > |\PPT_{W}(A)|.
\]
\item For any $W\subseteq\{1,\dots,l\}$, $C_{l}\le |\PPT_{W}(A)|\le C_{l+1}$.
\end{enumerate}
\end{cor}
\begin{proof}
Part~(1) is a direct consequence of Theorem~\ref{thm:(l+3)-gon}.
Part~(2) follows from part~(1) and the facts that
$|\PPT_{\emptyset}(A)|= C_{l}$ and $|\PPT_{\{1,\dots,l\}}(A)|=
C_{l+1}$ (the latter comes again from Theorem~\ref{thm:(l+3)-gon},
taking $W=A_I$).
\end{proof}
Note that $\PPT_\emptyset(A)$
 is in bijection to the set of triangulations of
the convex $(l+2)$-gon with vertices $A\setminus \{p\} =
\{0,1,\dots,l+1\}$, hence its cardinality is the Catalan number
$C_l$. Curiously enough, $\PPT_{A_I}(A)$ (that is, the set of
pointed pseudo-triangulations in which the tip $p$ is joined to
everything), has the cardinality of the next Catalan number
$C_{l+1}$. This follows from Theorem~\ref{thm:(l+3)-gon}, but was
first proved in Section~5.3 of~\cite{RSS03} (see also the remark and
picture on pp.~728--729). There, an associahedron with vertex set
$\PPT_{A_I}(A)$ is obtained, and it is regarded as a 1-dimensional
analog  of the construction of the polytope of pointed
pseudo-triangulations for a planar point set. Following this
analogy, we regard part (1) of Corollary~\ref{coro:singlechain-ppt}
as a 1-dimensional analog of Conjecture \ref{conj:monotone}.

\begin{cor}
\label{coro:singlechain-pt}
\begin{enumerate}
\item For every $W\subseteq\{1,\dots,l\}$, $ 2^{|W|}C_l\le|\PT_{W}(A)|\le 2^{|W|}C_{l+1}$.

\item $ 3^{l}C_l\le |\PT(A)| \le 3^{l}C_{l+1}. $ 

\item In
particular, $|\PT_{W}(A)|\in \Theta(2^{|W|}4^l\, l^{-\frac{3}{2}})$,
$|\PPT(A)|\in \Theta(8^l\, l^{-\frac{3}{2}})$ and $|\PT(A)|\in
\Theta(12^l\, l^{-\frac{3}{2}})$.
\end{enumerate}
\end{cor}

\begin{proof}
Part (1) comes from applying part (2) of Corollary~\ref{coro:singlechain-ppt} to each summand
in the expression $|\PT_{W}(A)|=\sum_{W'\subseteq W} |\PPT_{W'}(A)|$ of
Lemma~\ref{lemma:ppt}. Part (2) comes from adding the inequalities in part (1) for all the
subsets $W\subseteq\{1,\dots,l\}$. Finally, part (3) follows from $C_l\in \Theta(4^l
l^{-3/2})$.
\end{proof}
 Since $C_{l+1}/C_l < 4$, parts (1) and (2) of this corollary approximate the
numbers $|\PT(A)|$ and $|\PPT(A)|$ within a factor of four. In
Section~\ref{section:extras} we show how to obtain much better approximations.

\section{Proof of Theorem~\ref{thm:(l+3)-gon}}
\label{sec:conjecture3}

 Let us recall the statement we want to prove:

\bigskip\noindent
{\bf Theorem \ref{thm:(l+3)-gon}.}\ \ {\it
 Let $A$ be a single
chain with $l$ interior vertices, let $W\subseteq A_I$ be a subset of them.
 Let $B$ be the convex $(l+3)$-polygon with vertex set
$A\setminus\{p\} \cup \{q\}$, where $q$ is an extra point ``on the side opposite'' to
$p$.

Then, the pointed pseudo-triangulations of $A$ in which the interior neighbors of $p$ are
{\bf exactly} the points in $W$ (that is, the elements of $\PPT_W(A)$) are in bijection
to the triangulations of $B$ in which the interior neighbors of $q$ are {\bf contained}
in~$W$.} \medskip

\noindent Let us denote with $\T_W(B)$ the set of triangulations of
$B$ mentioned in the statement. The way we prove that $|\T_W(B)|$
and $|\PPT_W(A)|$ are the same number is by showing that both
families of numbers satisfy the same recursive formula. To this end,
let $W$ be a non-empty subset of $\{1,\dots,l\}$ and choose an element $v\in
W$. Let $W_1=\{ w\in W : w<v\}$ and $W_2=\{ w\in W : w>v\}$ be the
sets of elements of $W$ on both sides of $v$. Moreover, let:
\begin{itemize}
\item $A_1$ and $A_2$ be the ``single chains'' having as vertices $\{p,0,1,\dots,v\}$ and
$\{p,v,v+1,\dots,l+1\}$, respectively. \item $B_1$ and $B_2$ be the convex polygons
having as vertices $\{q,0,1,\dots,v\}$ and $\{q,v,v+1,\dots,l+1\}$, respectively.
\end{itemize}
Then:

\begin{lem}
\label{lemma:t-recursion} The following  recurrence holds:
\[
|\T_W(B)|-|\T_{W\setminus\{v\}}(B)|=|\T_{W_1}(B_1)|\cdot|\T_{W_2}(B_2)|.
\]
\end{lem}

\begin{proof}
The difference in the left-hand side coincides with the
triangulations of~$B$ that use the edge~$(q,v)$ and have the (other)
interior neighbors of~$q$ contained in~$W$. Clearly, those
triangulations can be obtained by triangulating $B_1$ and $B_2$
independently.
\end{proof}
 In the rest of this section we prove that the same recursion holds for the
numbers $|\PPT_{W}(A)|$, except we do it under the assumption that $v$ is the first
element in $W$. This assumption is enough for our purposes because knowing that
\[
|\PPT_W(B)|-|\PPT_{W\setminus\{v\}}(B)|=|\PPT_{W_1}(B_1)|\cdot|\PPT_{W_2}(B_2)|
\]
for any particular $v$, together with the inductive hypothesis that
$|\PPT_{W'}(A')|$ and $|\T_{W'}(B')|$ coincide whenever $|W'| <
|W|$ and the base case
$|\PPT_{\emptyset}(A)|=|\T_{\emptyset}(B)|$, implies that
$|\PPT_{W}(A)| =|\T_{W}(B)|$. That is, in order to prove
Theorem~\ref{thm:(l+3)-gon} it is enough to prove:

\begin{prop}
\label{prop:conjecture3reprise} For every $W\subseteq\{1,\dots,l\}$
and for $v=\min(W)$,
\[
|\PPT_W(A)|-|\PPT_{W\setminus\{v\}}(A)|=|\PPT_{W_1}(A_1)|\cdot |\PPT_{W_2}(A_2)|.
\]
\end{prop}

Our first observation is that:

\begin{lem}
\label{lemma:productPPTs} If $v=\min(W)$ then $|\PPT_{W_1}(A_1)|\cdot|\PPT_{W_2}(A_2)|$
equals the number of elements of $\PPT_{W}(A)$ that use the edge $(v,l+1)$.
\end{lem}

\begin{proof}
The edges $(p,v)$ and $(v,l+1)$ separate the triangle $\operatorname{conv}(A)$ into two
regions, so that we can count their number of pointed pseudo-triangulations
independently. The region on the right is the convex hull of $A_2$. Hence,
 it only
remains to show that the region on the left, let us denote it $A_L$, has the same number
of pointed pseudo-triangulations that join $W_1$ to $p$ as $A_1$ has. Note that
$v=\min(W)$ implies $W_1=\emptyset$. That is, $|\PPT_{W_1}(A_1)|$ is just the number of
triangulations of the $v+1$ points in convex position $\{0,\dots,v\}$. For $A_L$, we know
in addition that none of the vertices $\{0,\ldots,v-1\}$ can be connected to $v$, or $v$
would be non-pointed in the pseudo-triangulation of $A$ under consideration otherwise.
Thus, $|\PPT_{W_1}(A_L)|$ equals the number of triangulations of the $v+1$ points
$\{0,\ldots,v-1,l+1\}$, too.
\end{proof}

 Let $\PPT_{W}(A)^*$ denote the elements of $\PPT_{W}(A)$
that do not use the edge $(v,l+1)$. The above lemma implies that
Proposition~\ref{prop:conjecture3reprise} is equivalent to:

\begin{prop}
\label{prop:PPT*} For every $W\subseteq\{1,\ldots ,l\}$ and for
$v=\min(W)$,
\[
|\PPT_{W}(A)^*| = |\PPT_{W\setminus\{v\}}(A)|.
\]
\end{prop}
 We will prove this via an explicit (although complicated) bijection. For it, we
classify the elements of $\PPT_{W}(A)$ and $\PPT_{W\setminus\{v\}}(A)$ via the following
parameters.

\begin{defi}
\rm
Let $W=\{v_1,\dots,v_k\}$ and $v=\min(W)=v_1$.
\begin{enumerate}
\item For an element $T$ of $\PPT_{W}(A)$, we call \emph{end-point vector} of $T$ the vector
$(x_1,\dots,x_k)$ of length $|W|$ and with entries taken from
$\{0,\dots,l+1\}$, defined as follows: For every~$i$, there is a
single pseudo-edge in~$T$ having $v_i$ as a reflex vertex. Since
$v_i\in W$ and $T\in \PPT_W(A)$, one of the two corners joined by
this pseudo-edge is the tip~$p$. We define~$x_i$ to be the other
corner.
\item Similarly, the \emph{end-point vector} of an element $T$ of
$\PPT_{W\setminus\{v\}}(A)$ is the vector $(x_1,\dots,x_k)$ of length $|W|$ and with
entries in $\{0,\dots,l+1\}$ such that
\begin{itemize}
\item $x_1$ is the third corner of the triangle below the edge $(v-1,v)$. This is
well-defined because neither $v-1$ nor $v$ belong to $W\setminus\{v\}$, hence the edge
$(v-1,v)$ is in $T$ and both $v-1$ and $v$ are corners of the pseudo-triangle below it.
Moreover, this pseudo-triangle must necessarily be a triangle. \item For every $i>1$,
$x_i$ is defined as in part (1).
\end{itemize}
\end{enumerate}

\end{defi}

In the rest of this section we denote by $\PPT_{W,(x_1,\dots,x_k)}$ the
subset of $\PPT_{W}(A)$ consisting of elements with end-point vector $(x_1,\dots,x_k)$. Similarly, we denote
by $\PPT_{W\setminus\{v\},(x_1,\dots,x_k)}$ the elements of $\PPT_{W\setminus\{v\}}(A)$
with end-point vector $(x_1,\dots,x_k)$. With this notation, the main part of the proof
is to show the following bijections:

\begin{lem}
\label{lemma:bijections} Let $(x_1,\dots,x_k)$ be a vector of length $k$ with entries in
$\{0,\dots,l+1\}$.
\begin{enumerate}
\item If $x_1 > v$, then $
|\PPT_{W,(x_1,\dots,x_k)}|=|\PPT_{W\setminus\{v\},(x_1,\dots,x_k)}| $.

\item If $0<x_1 < v$, then $
|\PPT_{W,(x_1,\dots,x_k)}|=|\PPT_{W\setminus\{v\},(x_1-1,x_2^*,\dots, x_k^*)}| $, where
$x_i^*$ equals $x_i$ (respectively, equals $v$) if $x_i\ne x_1$ (resp., $x_i=x_1$).

\item If $x_1 =0$, then $ |\PPT_{W,(x_1,\dots,x_k)}|=|\PPT_{W\setminus\{v\},(l+1,
x_2^*,\dots, x_k^*)}| $, where $x_i^*$ equals $x_i$ (respectively, equals $v$) if $x_i\ne
x_1$ (resp., $x_i=x_1$).

\end{enumerate}
\end{lem}

\begin{proof}
(1) The bijection comes from a simple flip of the edge $(p,v)$ into
the edge $(v-1,v)$.

(2) As in part~(1), the first step is to perform a flip of the edge
$(p,v)$. This introduces an edge $(v,y)$, where $y=v+1$ unless both
$v_2=v+1$ and $x_2>v_2$ hold, in which case $y=x_2$. (See
Figure~\ref{fig:bijectionbeforerearrange}).

\begin{figure}[htb]
  \centering
  \includegraphics[width = \textwidth]{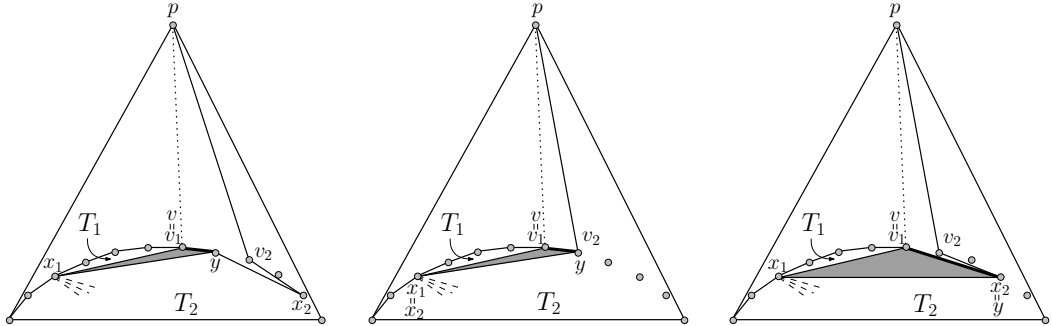}
  \caption{Three examples of flipping edge $(p,v)$ in part (2) of the proof of Lemma~\ref{lemma:bijections}.
  In the left and middle cases~$y=v+1$, in the right case~$y=x_2$.}
  \label{fig:bijectionbeforerearrange}
\end{figure}

In any case, we now have a pointed pseudo-triangulation $T$ that
belongs to $\PPT_{W\setminus\{v\}}$ and with the property that it
contains the triangle $t=(x_1,v,y)$. This triangle decomposes $T$
into three parts: a triangulation $T_1$ of the convex
$(v-x_1+1)$-gon with vertices $\{x_1,\dots,v\}$, the triangle
$(x_1,v,y)$ itself, and a pointed pseudo-triangulation $T_2$ of the
single chain with $l-(v-x_1)$ vertices $A\setminus
\{x_1+1,\dots,v\}$.

We are going to rearrange these three pieces in order to obtain a different pointed
pseudo-triangulation of $A$. We  embed $T_2$ as a pointed pseudo-triangulation of the
vertex set $A\setminus \{x_1,\dots,v-1\}$, add the triangle $(x_1-1,v-1,v)$ to it, and
then place the triangulation $T_1$ on the polygon $\{x_1-1,\dots,v-1\}$; see
Figure~\ref{fig:bijectionafterrearrange}. Essentially, in $T_2$ we are substituting
vertex $v$ for vertex $x_1$, and then we are changing the rest to be consistent with this
replacement. Since everything in $T_2$ previously joined to $x_1$ is now joined to $v$,
the new pointed pseudo-triangulation is indeed in
$\PPT_{W\setminus\{v\},(x_1-1,x_2^*,\dots,x_k^*)}$.

\begin{figure}[htb]
  \centering
  \includegraphics[width = \textwidth]{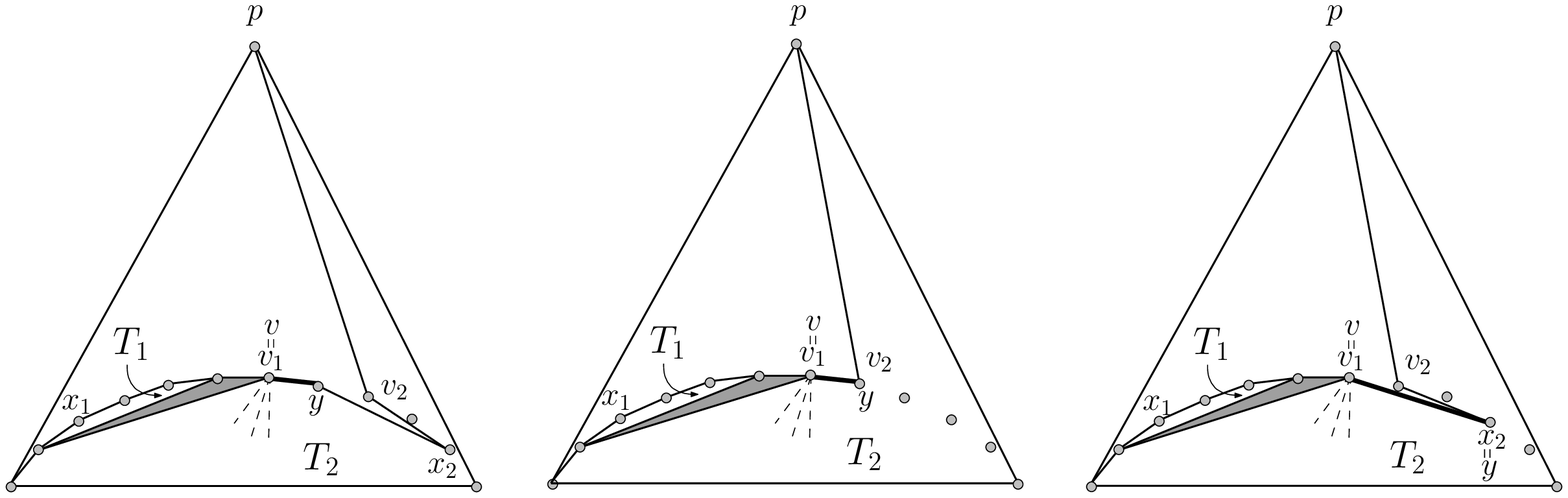}
  \caption{Rearrangement for the three examples in Figure~\ref{fig:bijectionbeforerearrange}.}
  \label{fig:bijectionafterrearrange}
\end{figure}

This process can be reversed: Starting with a pointed
pseudo-triangulation in \newline
$\PPT_{W\setminus\{v\},(y_1,y_2,\dots,y_k)}$, with
$y_1\in\{0,\dots,v-2\}$, the triangle $(y_1,v-1,v)$ decomposes it
into three parts: a triangulation $T'_1$ of the convex polygon with
vertices $\{y_1,\dots,v-1\}$, the triangle itself, and a pointed
pseudo-triangulation $T'_2$ of $A \setminus\{y_1+1,\dots,v-1\}$. We
place $T'_2$ on $A \setminus\{y_1+2,\dots,v\}$, $T'_1$ on
$\{y_1+1,\dots,v\}$ and insert the triangle $(y_1+1,v,v+1)$. We now
flip the edge $(v,v+1)$, and get a pointed pseudo-triangulation in
$\PPT_{W,(y_1+1,y_2^*,\dots,y_k^*)}$ where $y_i^*=y_i$ if $y_i\ne v$
and $y_i^*=y_1+1$ otherwise.

(3) The process is exactly the same as in part (2), except that
since $x_1=0$ we have to use $l+1$ in the role that was played by
$x_1-1$. That is, $T_1$ will be a triangulation of the convex
polygon $\{0,\dots, v\}$ before the rearrangement, and placed as a
triangulation $T'_1$ of the convex polygon $\{l+1,0,1,\dots, v-1\}$.
The rest is unchanged.
\end{proof}

\noindent
{\emph{Proof of Proposition~\ref{prop:PPT*}:\quad}}
We now show how Lemma~\ref{lemma:bijections} can be used to finish
the proof of Proposition~\ref{prop:PPT*}, hence that of
Proposition~\ref{prop:conjecture3reprise} (and therefore the one of
Theorem~\ref{thm:(l+3)-gon}). What we need to show is that the sets
in the right-hand sides of Lemma~\ref{lemma:bijections} cover the
set $\PPT_{W\setminus\{v\}}(A)$ without repetitions if we exclude
from the left-hand side the ones with $x_1=l+1$, which are the
elements in  $\PPT_{W}(A)\setminus \PPT_{W}(A)^*$. To this end, we
consider an end-point vector $Y=(y_1,\dots,y_k)$ of an element in
$\PPT_{W\setminus\{v\}}(A)$, and show that it comes from a unique
end-point vector $X=(x_1,\dots,x_k)$ of an element of
$\PPT_{W}(A)^*$ via the bijections in Lemma~\ref{lemma:bijections}.
First, observe that the end-points in the right-hand sides of parts
(1), (2) and (3) are distinguished by the properties $v<y_1<l+1$,
$y_1<v-1$ and $y_1=l+1$, respectively. Also, a valid $Y$ cannot have
$y_1$ equal to $v$ or $v-1$, by the definition of end-point vector
in $\PPT_{W\setminus\{v\}}(A)$. It only remains to show how to
recover the vector $X$ from $Y$:
\begin{enumerate}
\item[(a)] If $v < y_1 < l+1$, then just let $X=Y$.

\item[(b)] If $y_1<v-1$, then let $x_1=y_1+1$ and for $i>1$ let $x_i$ equal $x_1$ or
$x_i$ depending on whether $y_i=v$ or $y_i\ne v$.

\item[(c)] If $y_1 = l+1$, then let $x_1=0$ and let $x_i$ equal $x_1$ or $x_i$ depending
on whether $y_i=v$ or $y_i\ne v$.
\end{enumerate}
Note that, in all cases, $x_1\neq l+1$, trivially for (a) and (c)
and because $v-1\leq l-1$ in (b).
\qed
\bigskip

\section{Additional bounds and properties for the single chain}
\label{section:extras}

 Corollary~\ref{coro:singlechain-pt} gives the
approximations $|\PPT(A)|\simeq 2^lC_l$ and $|\PT(A)|\simeq 3^lC_l$
within a factor of four. In this section we show that
$|\PPT(A)|\simeq 2^{l+1}C_l$ and $|\PT(A)|\simeq 3^{l+1}C_l/2$ are
much better approximations, with errors of $12.5\%$ and $4\%$
respectively when~$l$ goes to infinity.
 This is in contrast with the fact
that we do not know such good and simple approximations for the individual summands
$|\PPT_W(A)|$. Our first step is to compute the sum of all the $\PPT_{W}(A)$'s for each
cardinality of~$W$, via the following recursive formulae.

\begin{thm}
\label{thm:singlechain} Let $ a(l,i):=\sum_{|W| = i} |\PPT_W(A)|. $ Then:
\begin{enumerate}
\item $a(l,0)=C_l$, and $a(l,1)=(l+1) C_l$. \item For every $i\ge 2$,
\[
a(l,i)={l+1\choose i} C_l -a(l-1,i-2).
\]
\end{enumerate}
\end{thm}

 As a preparation for the proof of
Theorem~\ref{thm:singlechain}, observe that the number ${l+1\choose
i} C_l$ that appears in the statement equals the number of ways of
specifying a triangulation of the $(l+2)$-gon together with~$i$ of
the $l+1$ boundary edges of the $(l+2)$-gon visible from the tip. We
say that a pointed pseudo-triangulation $T$ of $A$ is
\emph{compatible} with this specification if $T$ restricted to the
interior of the $(l+2)$-gon gives that triangulation and when
restricted to the boundary of the $(l+2)$-gon the~$i$ edges chosen
above are precisely the ones not appearing. Note that this notion of
compatibility is usable in both directions, i.e. ``pointed
pseudo-triangulations compatible with a choice'' and ``choices
compatible with a pointed pseudo-triangulation''. Some pointed
pseudo-triangulations of $A$ may not produce a triangulation of the
$(l+2)$-gon, and hence they are not compatible with any choice.
Reciprocally, some choices are not compatible with any pointed
pseudo-triangulation, but the next statement describes them:

\begin{lem}
\label{lemma:atmostone} Let a choice of a triangulation of the $(l+2)$-gon and a choice
of a subset of boundary edges of the $(l+2)$-gon visible from the tip be given. Then:
\begin{enumerate}
\item The choice is compatible with a pointed pseudo-triangulation of
  $A$ if and only if no ear of the triangulation is incident to two
  missing boundary edges.
\item A compatible choice determines uniquely a pointed
  pseudo-triangulation.  This pseudo-triangulation uses $i$ interior
  edges incident to the tip vertex, where $i$ equals the number of
  missing boundary edges.
\end{enumerate}
\end{lem}

\begin{proof}
  Observe that compatible means that from the given choice of triangulation and
  subset of boundary edges we can get a pointed
  pseudo-triangu\-la\-tion, by adding to the chosen triangulation some edges incident to the tip
  and removing the chosen boundary edges.  We
  call an ear of the triangulation incident to two missing boundary
  edges of the $(l+2)$-gon a {\em bad ear} (see
  Figure~\ref{fig:badear-pt}).

  \begin{figure}[htb]
    \centering
    \includegraphics{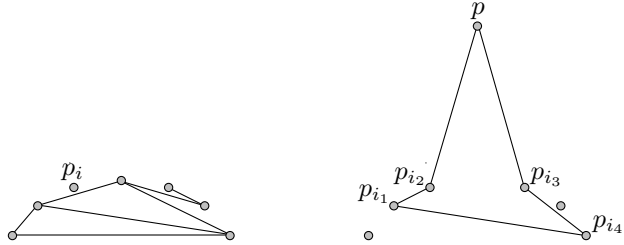}
    \caption{Left: bad ear at $p_i$. Right: bad pseudo-triangle in the proof of Theorem~\ref{thm:singlechain}.}
    \label{fig:badear-pt}
  \end{figure}
  Clearly, if a bad ear appears at $p_i$ the choice cannot be
  compatible with a pointed pseudo-triangulation, because a vertex
  cannot have degree 1 in a pointed pseudo-triangulation. Hence, assume
  that we have a choice with no bad ears and let us prove that it is
  compatible with one and only one pointed pseudo-triangulation. The
  way to obtain the pointed pseudo-triangulation is as follows: let
  $p_i p_{i+1}$ be a missing edge in the choice.  Let $p_k$ be the
  vertex of the triangulation of the $(l+2)$-gon joined to it. We add the
  edge $(p,p_i)$ or $(p,p_{i+1})$ depending on whether $k<i$ or $k>i+1$. The
  assumption of no bad ears implies that we add as many edges as
  missing edges were in the choice.  In particular, the set of edges
  obtained in this way has cardinality $2l+3$, the same as an element of
  $\PPT_\emptyset(A)$.  Since every pointed graph with $2l+3$ edges of
  a vertex set of size $l+3$ is a pointed pseudo-triangulation, we
  have shown existence.  For uniqueness, just observe that
  every compatible pointed pseudo-triangulation must have at least the
  edges we have added: an interior vertex not joined to the tip must
  be joined to vertices both to its right and to its left.
\end{proof}

\noindent
{\emph{Proof of  Theorem~\ref{thm:singlechain}:\quad}}
  The equation $a(l,0)=|\PPT_\emptyset(A)| =C_l$ is obvious. For $a(l,1)$,
  observe that every pointed pseudo-triangulation with a single
  interior edge joined to the tip uniquely gives rise, by a flip of that edge,
  to another one with no
  edges joined to the tip. Conversely, every pointed
  pseudo-triangulation with no edges to the tip gives rise to $l+1$
  pointed pseudo-triangulations with a single edge to the tip, by the
  $l+1$ possible flips of the boundary edges of the $(l+2)$-gon.

  For the proof of part (2), let us call $b(l,i):={l+1\choose i} C_l$ and let $c(l,i)$ be
  the number of pointed pseudo-triangulations of $A$ that are
  compatible with a choice of triangulation and boundary. We will
  abuse notation and use $a(l,i)$, $b(l,i)$ and $c(l,i)$ to represent
  not only the numbers but also the sets of objects counted by them; $\cup_{|W| = i} \PPT_W(A)$ for
  $a(l,i)$, choices of a triangulation and a subset of edges as above for $b(l,i)$, and the double meaning of ``pointed
  pseudo-triangulations compatible with a choice'' and ``choices
  compatible with a pointed pseudo-triangulation'' for~$c(l,i)$.

  Clearly, every element of $b(l,i)$ with $k$ ``bad ears" can be
  considered a member of $c(l-k,i-2k)$: just delete the $k$ interior
  points where the bad ears occur. Reciprocally, each member of~$c(l-k,i-2k)$
  can give a member of $b(l,i)$ in ${l-k-(i-2k)+1
  \choose k}$ ways: we choose $k$ of the $l-k-(i-2k)+1$ used
  boundary edges and place a new vertex (a bad ear) beyond each of
  those $k$ edges.  Hence:
\[
b(l,i) = \sum_{k\ge 0} {l-i+1+k \choose k} c(l-k,i-2k).
\]

Now, what can make a pointed pseudo-triangulation not compatible
with a choice of triangulation plus boundary edges is the existence
of a \emph{bad pseudo-triangle}
$[p,p_{i_2},p_{i_1},p_{i_4},p_{i_3}]$ with $i_1<i_2<i_3<i_4$, see
Figure~\ref{fig:badear-pt} (note that $p_{i_3} = p_{i_2}+1$ in order
to be a pseudo-triangle). In this case the restriction to the
$(l+2)$-gon has a quadrangle, let us call it a \emph{bad quadrangle}
instead of being a triangulation. More that one bad quadrangle can
occur, but the two edges $(p,p_{i_2})$ and $(p,p_{i_3})$ that join
one bad quadrangle to the tip cannot join any other bad quadrangle
to the tip.

In particular, if an element of $a(l,i)$ produces $k$ bad
quadrangles, contracting the edge $(p_{i_2},p_{i_3})$ of each
quadrangle and removing the $2k$ corresponding edges incident to the
tip we get an element of $c(l-k,i-2k)$, because of
Lemma~\ref{lemma:atmostone}. To get back an element of $a(l,i)$ from
one of $c(l-k,i-2k)$ one must choose $k$ of the $l-i+k$ interior
vertices not incident to the tip and split them into two vertices,
joining both to the tip. Clearly, there are ${l-i+k \choose k}$ ways
to do that. Hence:
\[
a(l,i) = \sum_{k\ge 0} {l-i+k \choose k} c(l-k,i-2k).
\]
Then:
\[
a(l-1,i-2) = \sum_{k\ge 1} {l-i+k \choose k-1} c(l-k,i-2k),
\]
where the index~$k$ has been shifted by one after evaluating with the previous formula.
To get the statement, add the two last equalities and compare them to the one for
$b(l,i)$.
\qed
\bigskip

\begin{table}[htb]
\begin{center}
\small
\begin{tabular}{|r|rrrrrr|c|}
\hline
 ${}_l \setminus {}^i$
        & 0 &1 &2 &3 &4 &5 &   $ |\PPT(A)|=\sum a(l,i)$\\
\hline
0&       1 &        &       &         &        &      &       1\\
1&       1 &     2 &       &         &        &      &       3\\
2&       2 &     6 &    5 &         &        &      &       13\\
3&       5 &   20 &  28 &    14 &        &      &       67\\
4&     14 &   70 & 135 & 120 &   42 &      &      381\\
5&      42 & 252 & 616 & 770 & 495&132 &      2307\\
\hline
\end{tabular}
\caption{Values of $a(l,i)$ and $|\PPT(A)|$ for $l,i \leq 5$.} \label{tab:ali}
\end{center}
\end{table}

 Theorem~\ref{thm:singlechain} allows us to compute all the
values of $a(l,i)$ recursively, starting from those stated in part
(1). The first few values of $a(l,i)$ are shown in
Table~\ref{tab:ali}. The recursion also tells us that the array
$a(l,i)$ equals the sequence A062991 in Sloane's Encyclopedia
\cite{integersequences}.  The row sums, that is, the numbers
$|\PT_{A_I}(A)|=|\PPT(A)|$ of all pointed pseudo-triangulations,
form the sequence A062992 and satisfy:
\[
|\PPT(A)| = |\PT_{A_I}(A)| = \sum_{i=0}^l a(l,i) = 2 \sum_{j=0}^l (-1)^{l-j} C_j 2^{j} -
(-1)^l .
\]
We can obtain them by adding over all values of $i$ in the formula
of Theorem~\ref{thm:singlechain}.

\begin{cor}
\label{coro:single} The number $a_l =|\PPT(A)|$ of pointed pseudo-triangulations of the
single chain satisfies:
\[
a_l = 2^{l+1} C_l - a_{l-1}.
\]
Hence,
\[
\left(1-\sum_{i=\lfloor 2+\frac{(-1)^l}{2}\rfloor}^{l}(-1)^{l-i}\prod_{j=i}^{l}\frac{j+1}{4(2j-1)}\right)
\cdot 2^{l+1}C_l  \le |\PPT(A)| \le 2^{l+1}C_l.
\]
Observe that the parenthesis in the left-hand side tends to~$\frac{8}{9}$ when~$l$ goes to infinity.
\end{cor}

\begin{proof}
  The first statement follows from $|\PPT(A)| = \sum_{i=0}^l
  a(l,i)$ and Theorem~\ref{thm:singlechain}, using that $a(l,l)=C_{l+1}$
  (Theorem~\ref{thm:(l+3)-gon}, with $W=\{1,\dots,l\})$.
 For an example, $381 = 2^5 C_4 - 67 = 32\cdot 14 -67$. For the second part, the
 upper bound is straightforward and for the lower bound the first part gives
 $$a_l=(-1)^l+\sum_{i=1}^l(-1)^{l+i}2^{i+1}C_i$$
 and then one can use the fact that $C_{l}=C_{l-1}(4l-2)/(l+1)$.
\end{proof}

 We now turn our attention to the total number of pseudo-triangulations
 $\PT(A)$. Lemma~\ref{lemma:ppt} implies that:
\[
|\PT(A)| = \sum_{W'\subseteq A_I} 2^{|A_I\backslash
W'|}|\PPT_{W'}(A)| = \sum_{i=0}^l 2^{l-i} a(l,i).
\]

\begin{cor}
\label{coro:singlept} The number $b_l = |\PT(A)|$ of pseudo-triangulations of the single
chain satisfies:
\[
2 b_l = 3^{l+1} C_l - b_{l-1}.
\]
Hence,
\[
\left(1-\sum_{i=\lfloor 2+\frac{(-1)^l}{2}\rfloor}^{l}(-1)^{l-i}\prod_{j=i}^{l}\frac{j+1}{12(2j-1)}\right)
\cdot \frac{3^{l+1}}{2}C_l  \le |\PPT(A)| \le \frac{3^{l+1}}{2}C_l.
\]
Observe that the parenthesis in the left-hand side tends to~$\frac{24}{25}$ when~$l$ goes to infinity.
\end{cor}

\begin{proof}
  Similar to the proof of Corollary~\ref{coro:single}. Fo the second part we use
  $|\PT(A)| = \sum_{i=0}^l 2^{l-i} a(l,i)$ and the first part to get
  $$b_l=\frac{1}{(-2)^l}+\sum_{i=1}^l(-1)^{l+i}\frac{3^{i+1}}{2^{l-i+1}}C_i.$$
\end{proof}

\section{The double chain}
\label{sec:doublechain}

\begin{figure}[htb]
  \centering
  \includegraphics{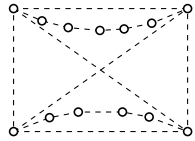}
  \caption{ A double chain: $l = 5$ and $m = 4$.}
  \label{fig:double}
\end{figure}

 For any two numbers $l,m\ge 0$, we call {\em double chain} with parameters
$(l,m)$ the point set consisting of a convex 4-gon with $l$ and $m$ points, respectively,
placed forming concave chains next to opposite edges of the 4-gon in a way that they do
not cross the two diagonals of the convex 4-gon (see Figure~\ref{fig:double}). The double
chain decomposes into a convex $(l+2)$-gon, a convex $(m+2)$-gon, and a non-convex
$(l+m+4)$-gon, the latter  with ${l+m+2 \choose l+1}$ triangulations~\cite{GNT00}. Hence,
the double chain has exactly
\[
C_{l}C_{m}{l+m+2 \choose l+1}
\]
triangulations. In the extremal case $l=m=(n-4)/2$ this gives $\Theta(8^{n} n^{-7/2})$.
The double chain has been, until very recently~(see~\cite{ahhhkv06}), the example of a
point set in the plane with asymptotically the biggest number of triangulations known.

Throughout this section, let $A$ be a double chain with $l$ and $m$ interior points in
the two chains, respectively (so $A$ has $l+m+4$ points in total). We call the $l+2$ and
$m+2$ vertices in the two chains the ``top'' and ``bottom'' parts.

\begin{figure}[htb]
  \centering
  \subfigure[$T_A$]{\includegraphics{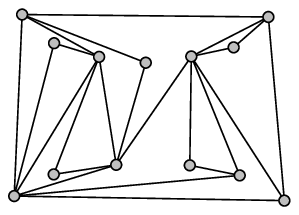}}
  \hspace{1in}
  \subfigure[$T_B$ and $T_C$]{\includegraphics{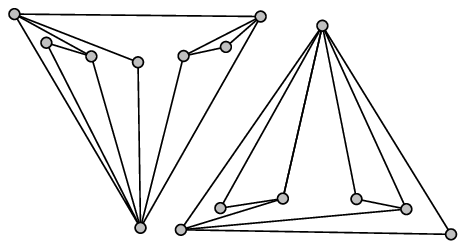}}
  \caption{Decomposing a pseudo-triangulation of a double chain.}
  \label{fig:doublepseudo}
\end{figure}

In order to count the number of pseudo-triangulations of $A$, let us
call $B$ and $C$ single chains with $l$ and $m$ interior points
each. $B$ can be considered the subset of $A$ consisting of the top
part plus a bottom vertex, and analogously for $C$. Every
pseudo-triangulation $T_A$ of $A$ induces on the one hand a
pseudo-triangulation~$T_B$ of~$B$ by contracting all bottom vertices
to a single one, and on the other hand a pseudo-triangulation~$T_C$
of~$c$ by doing the same with all top vertices (see
Figure~\ref{fig:doublepseudo}). Since no pseudo-triangle of~$T_A$
contains both more than one top vertex and more than one bottom
vertex, every pseudo-triangle survives either in~$T_B$ or in~$T_C$
but not in both.

Conversely, given a pair of pseudo-triangulations $T_B$ and $T_C$ of $B$ and $C$, if $i$
(resp.  $j$) denotes the number of interior edges incident to the bottom point in $T_B$
(resp.  to the top point in $T_C$), there are exactly ${i+j+2\choose i+1}$ ways to
recover a pseudo-triangulation of $A$ from that data, by shuffling the $i+1$
pseudo-triangles of $T_B$ incident to the bottom and the $j+1$ of $T_C$ incident to the
top.

\begin{thm}
\label{thm:doublechain} Let $V$ and $W$ be subsets of the top and bottom interior points.
For each $i\le v\le l$ and $j\le w\le m$ let $t^{v,w}_{i,j} :={l-
  v+i + m- w+j +2 \choose l-
  v+i+1}$.  Then:
\begin{enumerate}
\item
\[
|\PT_{V\cup W} (A)| = \sum_{{}^{V'\subseteq V}_{W'\subseteq W}}
  t^{|V|,|W|}_{|V'|,|W'|} |\PPT_{V'}(B)| |\PPT_{W'}(C)|.
\]
\item In particular, if $v=|V|$ and $w=|W|$ then,
\[
|\PT_{V\cup W} (A)| \in \Theta\left(C_l C_m \sum_{i=0}^v \sum_{j=0}^w {v\choose
i}{w\choose j}  t^{v,w}_{i,j}\right).
\]
\end{enumerate}
\end{thm}

\begin{proof}
  The first observation is that the ``shuffling'' described above
  preserves pointedness.  Then, part (1) follows from the fact that
 in the expression
\[
|\PT_V(B)| =\sum_{V'\subseteq V} |\PPT_{V'}(B)|
\]
of Lemma \ref{lemma:ppt}, each element of $\PPT_{V'}(B)$ corresponds to an element of
$\PT_V(B)$ with exactly $l-|V\backslash V'|=l-|V|+|V'|$ interior edges incident to the
bottom point (same for $C$).

Part (2) follows from part (1) using that $|\PPT_{V'}(B)|\in \Theta(C_l)$,
$|\PPT_{W'}(C)|\in \Theta(C_m)$, and in the sum of part (1) there are exactly ${v\choose
i}{w\choose j}$ summands with $|V'|=i$ and $|W'|=j$.
\end{proof}

\begin{cor}\label{cor:conjecturedouble}
The double chain satisfies Conjecture~\ref{conj:monotone}.
\end{cor}

\begin{proof}
 When we add a point $p$ to $V$, Theorem~\ref{thm:doublechain} gives that
$|\PT_{V\cup\{p\}\cup W}(A)|$ equals
  \[
 \sum_{{}^{V'\subseteq V}_{W'\subseteq W}}
\left( t^{|V\cup\{p\}|,|W|}_{|V'|,|W'|}
  |\PPT_{V'}(B)|
  +
  t^{|V\cup\{p\}|,|W|}_{|V'\cup\{p\}|,|W'|} |\PPT_{V'\cup\{p\}}(B)|
\right)
  |\PPT_{W'}(C)|.
\]
We neglect the first summand, and use monotonicity
 of $|\PPT_{V'}|$ (part (1) of
Corollary~\ref{coro:singlechain-ppt}) in the second summand, which is then greater than
  \[
\sum_{{}^{V'\subseteq V}_{W'\subseteq W}}
  t^{|V\cup\{p\}|,|W|}_{|V'\cup\{p\}|,|W'|} |\PPT_{V'}(B)|
  |\PPT_{W'}(C)|.
\]
But this equals $|\PT_{V\cup  W}(A)|$ by Theorem~\ref{thm:doublechain} since, clearly,
  \[
  t^{v,w}_{i,j}=
  t^{v+1,w}_{i+1,j}.  \]
\end{proof}

 Part (2) of Theorem~\ref{thm:doublechain} implies that, to understand the
asymptotics of pseudo-triangulations of the double chain, we need to understand the
expressions
\begin{equation}\label{eqn:doublechain}
\sum_{i=0}^v \sum_{j=0}^w {v\choose i}{w\choose j}  t^{v,w}_{i,j}=
\sum_{i=0}^v \sum_{j=0}^w {v\choose i}{w\choose j}  {l-i+m-j+2
\choose l-i+1}.
\end{equation}
The second form is obtained from the first by the substitutions $i\to v-i$ and $j\to
w-j$, since $t^{v,w}_{v-i,w-j} = {l-i+m-j+2 \choose l-i+1}$.

For the special case $l=v$ and $m=w$, this expression has a very nice combinatorial
interpretation and has appeared in the literature (see below). In particular, we can give
the exact asymptotics of the number of pointed pseudo-triangulations of a double chain
with $l=m$ (Theorem~\ref{thm:doublechain-coker}). For general values of $l$, $m$, $|V|$,
$|W|$, or for the total number of pseudo-triangulations, we can only offer the upper and
lower bounds in the following two statements:

\begin{lem}
\label{lemma:doublechain-ub}
\begin{enumerate}
\item For every $V\subseteq\{1,\dots, l\}$ and $W\subseteq \{1,\dots,m\}$,
\[
|\PT_{V\cup W} (A)| \in O\left(2^{l+m}\, \left({3/2}\right)^{v+w}\,C_{l} C_{m}\right),
\]
where $v=|V|$ and $w=|W|$. \item In particular,
\[
|\PPT(A)| \in O(3^{l+m}\,C_{l} C_{m}) =  O(12^{l+m}\,(lm)^{-3/2}).
\]
\item
\[
|\PT(A)| \in O(5^{l+m}\,C_{l} C_{m}) =  O(20^{l+m}\,(lm)^{-3/2}).
\]
\end{enumerate}
\end{lem}

\begin{proof}
Starting with the equality in Theorem~\ref{thm:doublechain}, we bound $ |\PPT_{V'}(B)| $
by $C_{l+1}$and $|\PPT_{W'}(C)|$ by $C_{m+1}$, respectively, using part (2) of
Corollary~\ref{coro:singlechain-ppt}. We also bound
\[
 t^{v,w}_{|V'|,|W'|} = {l- v+| V'| + m- w+| W'| +2 \choose l-  v+|V'|+1}
\le 2^{l+m-v-w +  |V'| + |W'| +2}.
\]
Thus,
\begin{equation}\label{eqn:asymptotics_PT_double_chain}
\frac{|\PT_{V\cup W} (A)|}{C_{l+1} C_{m+1}} \le  2^{{l+m-v-w +2}}
\sum_{{}^{V'\subseteq V}_{W'\subseteq W}} 2^{|V'|}\,2^{|W'|} =
2^{{l+m-v-w +2}}  3^{v}\,3^{w}.
\end{equation}
That finishes part (1), since $C_{l+1}\in \Theta(C_l)=\Theta(4^l
l^{-3/2})$. For the upper bound in part (2), we simply specialize
$v=l$ and $w=m$. For the upper bound in part (3) we add over all
values of $V$ and $W$ the
inequality~(\ref{eqn:asymptotics_PT_double_chain}) obtained above,
since
\[
|\PT(A)| =
\sum_{{}^{V\subseteq\{1,\dots,l\}}_{W\subseteq\{1,\dots,m\}}}
|\PT_{V\cup W} (A)|.
\]
Hence:
\[
\frac{|\PT (A)|}{C_{l+1} C_{m+1}} \le \sum_{v=0}^{l}\sum_{w=0}^{m} {l \choose v}{m
\choose w }2^{{l+m-v-w +2}}  3^{v+w} =
\]
\[
=2^{{l+m+2}} \sum_{v=0}^{l}\sum_{w=0}^{m} {l \choose v}{m \choose w }
    \left(\frac{3}{2}\right)^{v+w} =
    2^{{l+m+2}}  \left(\frac{5}{2}\right)^{l+m} =
    4\cdot{5}^{l+m}.
\]
\end{proof}

We now look at lower bounds. We obtain the following ones by simply taking the greatest
summand in the expressions derived from Theorem~\ref{thm:doublechain}. Observe that in
the case $l=m$ they differ from the upper bounds only by a polynomial factor of
$l^{-3/2}$ and $l^{-5/2}$, respectively.

\begin{thm}
\label{thm:doublechain-lb}
\begin{enumerate}
\item $\displaystyle
 |\PPT(A)| \in
            \Omega\left(3^{l+m} C_l C_m
             \frac{(l+m)^{1/2}}{lm}
             \left(\frac{1}{2}\right)^{2|l-m|/3}
              \right).
              \medskip
$ \item $\displaystyle
 |\PT(A)| \in
            \Omega\left(5^{l+m} C_l C_m
             \frac{(l+m)^{1/2}}{(lm)^{3/2}}
             \left(\frac{1}{2}\right)^{4|l-m|/5}
              \right).
$
\medskip
\item In particular, if $l=m=(n-4)/2$ (where $n$ is the total number of vertices), we
have
\[
|\PPT(A)| \in \Theta^*(12^n), \qquad \text{and} \qquad |\PT(A)| \in
\Theta^*(20^n).
\]
\end{enumerate}
\end{thm}

\begin{proof}
For part~(1) we start with
\[
|\PPT(A)|=\sum_{{}^{V'\subseteq \{1,\dots,l\}}_{W'\subseteq \{1,\dots,m\}}}
t^{l,m}_{|V'|,|W'|}\ |\PPT_{V'}(B)|\ |\PPT_{W'}(C)| \ge
\]
\[
\ge \sum_{{}^{V'\subseteq \{1,\dots,l\}}_{W'\subseteq \{1,\dots,m\}}}
t^{l,m}_{|V'|,|W'|}\ C_l\, C_m =C_l\, C_m \sum_{i=0}^l \sum_{j=0}^m {l\choose i}{m\choose
j} {i+j+2 \choose i+1}.
\]
In this expression we substitute the sum by the summand with $i=2l/3$ and $j=2m/3$. That
is:
\[
\frac{|\PPT(A)|}{C_l C_m}
{l\choose 2l/3}{m\choose 2m/3} {\frac{2(l+m)}{3}+2 \choose 2l/3} \sim
{l\choose 2l/3}{m\choose 2m/3} {2(l+m)/3 \choose 2l/3}.
\]
Next we approximate the binomial coefficients using Stirling approximation, which gives:
\[
{l\choose 2l/3}\in \Theta\left(\frac{3^l}{2^{2l/3}}l^{-1/2}\right)
\]
and
\[{2(l+m)/3 \choose 2l/3}\in \Theta\left(
\left(\frac{l+m}{l}\right)^{2l/3}
             \left(\frac{l+m}{m}\right)^{2m/3}
             \left(\frac{l+m}{lm}\right)^{1/2}
\right).
\]
Putting things together we get

\[
\frac{|\PPT(A)|}{C_lC_m} \in \Omega\left(3^{l+m}
             \left(\frac{l+m}{2l}\right)^{2l/3}
             \left(\frac{l+m}{2m}\right)^{2m/3}
             \frac{(l+m)^{1/2}}{lm}
              \right)
 =
 \]
 \[=
 \Omega\left(3^{l+m}
              \left(\frac{(l+m)^2}{4lm}\right)^{\frac{2\min(l,m)}{3}}
             \left(\frac{l+m}{2\max(l,m)}\right)^{\frac{2|l-m|}{3}}
             \frac{(l+m)^{1/2}}{lm}
 \right).
\]
This gives part~(1), since $\frac{(l+m)^2}{4lm}\ge 1$ and
$\frac{l+m}{2\max(l,m)}\ge \frac{1}{2}$.

For part (2) we use the same ideas. We start with
\[
|\PT(A)| \ge 
\sum_{v=0}^l\sum_{w=0}^m {l\choose v}{m\choose w}
    \sum_{i=0}^v\sum_{j=0}^w {v\choose i}{w\choose j}
                          {\scriptstyle{l-v+i+m-w+j+2} \choose \scriptstyle{l-v+i+1}} C_l C_m.
\]
Here, we substitute the sum with the summand $i=2l/5$, $v=3l/5$, $j=2m/5$, and $w=3m/5$.
This gives:
\[
\frac{|\PT(A)|}{C_lC_m} \ge {l\choose 3l/5}{m\choose 3m/5}
             {3l/5\choose 2l/5}{3m/5\choose 2m/5}
                          {4(l+m)/5\choose 4l/5} .
\]
As before, Stirling's approximation gives:
\[
{l\choose 3l/5} \in\Theta\left(\frac{5^{l}}{3^{3l/5}2^{2l/5}}l^{-1/2}\right),
\qquad\qquad {3l/5\choose 2l/5} \in\Theta\left(\frac{3^{3l/5}}{2^{2l/5}}l^{-1/2}\right)
\]
\[
{4(l+m)/5\choose 4l/5} \in \Theta\left( \left(\frac{l+m}{l}\right)^{4l/5}
             \left(\frac{l+m}{m}\right)^{4m/5}
             \left(\frac{l+m}{lm}\right)^{1/2}
\right).
\]
That is,
\[
\frac{|\PT(A)|}{C_lC_m}
    \in \Omega\left(5^{l+m}
             \left(\frac{l+m}{2l}\right)^{4l/5}
             \left(\frac{l+m}{2m}\right)^{4m/5}
             \frac{(l+m)^{1/2}}{(lm)^{3/2}}
 \right)=
\]
\[=
 \Omega\left(5^{l+m}
              \left(\frac{(l+m)^2}{4lm}\right)^{\frac{4\min(l,m)}{5}}
             \left(\frac{l+m}{2\max(l,m)}\right)^{\frac{4|l-m|}{5}}
             \frac{(l+m)^{1/2}}{(lm)^{3/2}}
 \right).
\]

Part~(3) is straightforward from parts~(1) and~(2), by
Lemma~\ref{lemma:doublechain-ub}.
\end{proof}

 Let us now restrict our attention to the case of pointed pseudo-triangulations.
Applying Theorem~\ref{thm:doublechain} with $v=l$ and $w=m$ we get the following, which
has been used in the proofs of the last two results:
\[
\frac{|\PPT(A)|}{C_lC_m}\in \Theta\left(\sum_{i=0}^l \sum_{j=0}^m {l\choose i}{m\choose
j}  {i+j+2 \choose i+1} \right).
\]
Let us call $E^{l,m}$ the expression inside the $\Theta(-)$. It
turns out that $E^{l,m}$ has the following nice interpretation: it
equals the number of lattice paths from $(0,0)$ to $(l+1,m+1)$ when
horizontal and vertical steps of arbitrary positive length are
allowed. In other words, it equals the number of monotone rook paths
from $(0,0)$ to $(l+1,m+1)$ (a path is specified not only by the
squares traversed, but also by the positions where the rook stops.
The rook is allowed to do several consecutive horizontal or vertical
moves). Indeed, for a particular path, $i+1$ and $j+1$ represent the
numbers of horizontal and vertical moves taken by the rook. The
coefficient ${l\choose i}$ (resp. ${m\choose j}$) accounts for the
possibilities of columns (resp. rows) where the rook makes at least
one stop, and the coefficient ${i+j+2 \choose i+1}$ accounts for the
relative ordering of the $i+1$ horizontal and $j+1$ vertical moves.

The sequence $E^{l,m}$ appears (with a shift in the indices) as A035002 in
\cite{integersequences} and has been studied in Section 7 of \cite{coker}. It satisfies,
among others, the following formulas:
\[
E^{l+1,m+1}=2E^{l,m+1}+2E^{l+1,m}-3E^{l,m}\quad \forall l,m>1; \quad
E^{l,0}=E^{0,l}=(l+4)2^{l-1},
\]
or
\[
\sum_{l+m=n} E^{l,m} = 2(3^{n+1}-2^{n+1}).
\]

In particular, the generating function of its diagonal sequence $E^{m,m}$ (sequence
A051708) is known, and from it we can derive the asymptotics very precisely:

\begin{thm}
\label{thm:doublechain-coker} Let $A$ be a double chain with $n$
vertices in total and with $l=m=(n-4)/2$. Then, $|\PPT(A)| \in
\Theta(12^n n^{-7/2})$.
\end{thm}

\begin{proof}
The generating function of $E^{m,m}$ is
\[
f(t)= \frac{9t-1 + \sqrt{9t^2-10t+1}}{2(9t-1)}= \frac{1}{2} + \frac{1}{2}
\sqrt{\frac{1-t}{1-9t}}.
\]
(This is Theorem 7.1(c) of \cite{coker}, except there a negative sign is wrongfully taken
before the square root. The correct sign is positive since otherwise $f(t)$ is negative
near zero, which does not make sense).

The dominant (i.e., smallest in absolute value) singularity of $f(t)$ is at $t=1/9$, and
near the singularity one has
\[
f(t)\sim \sqrt{\frac{2}{9}}\,(1-9t)^{-1/2}.
\]
Then, the singularity analysis of \cite{fo-90} (see also~\cite{fs-07}) implies that
\[
E^{m,m}\sim \sqrt{\frac{2}{9}} \cdot \frac{9^m}{\sqrt{\pi m}}.
\]
Hence, $|\PPT(A)| \in \Theta(\frac{9^m
{C_m}^2}{\sqrt{m}})=\Theta(12^n n^{-7/2})$.
\end{proof}

 A similar analysis could be undertaken for $\PT(A)$. By
Theorem~\ref{thm:doublechain} and the
equality~(\ref{eqn:doublechain}) we have $|\PT(A)|/C_lC_m\in \Theta
(F^{m,n})$, where
\[
F^{m,n}=\sum_{^{\,\, 0\le i\le v\le l}_{0\le j\le w\le m}} {l\choose v}{m\choose
w}{v\choose i}{w\choose j}  {l+m-i-j+2 \choose l-i+1}.
\]

$F^{m,n}$  can still be interpreted (although less directly) in terms of rook paths, and
satisfies formulas such as
\[
F^{l+1,m+1}=3F^{l,m+1}+3F^{l+1,m}-5F^{l,m}\ \forall l,m>1, \qquad
F^{l,0}=F^{0,l}=(l+6)3^{l-1},
\]
or
\[
\sum_{l+m=n} F^{l,m} = 5^{n-1}-3^{n-1}.
\]
The latter, since the biggest summand is obtained with $l=m=n/2$, implies that $F^{n,n}$
is between $\Omega(5^n n^{-1})$ and $O(5^n)$, in agreement with---but also refining---the
result in part (3) of Theorem~\ref{thm:doublechain-lb}. We believe that $F^{n/2,n/2}\in
\Theta(5^n n^{-1/2})$ and, hence, that the total number of pseudo-triangulations of a
double chain with the same number of points on both sides is in $\Theta(20^n n^{-7/2})$.


\bibliographystyle{abbrv}

\end{document}